\documentclass[leqno,11pt]{article}

\usepackage{amssymb}
\usepackage{amsthm}
\usepackage{eucal}
\usepackage{amsmath,amsxtra}

\makeatletter

\skip\footins 24pt plus 4pt minus 2pt

\makeatother

\newcommand{\fa}{{\mathfrak a}}
\newcommand{\fg}{{\mathfrak g}}
\newcommand{\fk}{{\mathfrak k}}
\newcommand{\fn}{{\mathfrak n}}

\newcommand\A{\mathcal A}

\newcommand\D{\mathcal D}
\newcommand\E{\mathcal E}
\newcommand\F{\mathcal F}
\renewcommand\P{\mathcal P}

\newcommand\s{\mathcal S}

\newcommand{\CC}{\mathbf{C}}

\newcommand{\DD}{\mathbf{D}}
\newcommand{\GG}{\mathbf{G}}
\newcommand{\HH}{\mathbf{H}}
\newcommand{\MM}{\mathbf{M}}
\newcommand{\OO}{\mathbf{O}}
\newcommand{\RR}{\mathbf{R}}

\renewcommand{\SS}{\mathbf{S}}
\newcommand{\ZZ}{\mathbf{Z}}

\newcommand{\ad}{\mathop{\rm ad\,}\nolimits}
\newcommand\Cinf{\mathcal{C}^\infty}

\newcommand{\Hom}{\mathop{\rm Hom\,}\nolimits}

\newcommand{\Tr}{\mathop{\rm Tr\,}\nolimits}
\newcommand{\Trace}{\mathop{\rm Trace\,}\nolimits}

\def \tilde{\widetilde}
\def \hat{\widehat}

\newcommand\bigcheck[1]{#1 \raise1ex\hbox{$\hspace{-1ex}{}^\vee$}}
\newcommand\sucheck[1]{#1 \raise0.5ex\hbox{$\hspace{-1ex}{}^\vee$}}

\newcommand{\checkvarphi}{\smash{\raisebox{-0.8ex}{\kern0.4ex
    \LARGE$\check{\smash{\raisebox{0.5ex}{\kern-0.2ex\normalsize
          $\varphi$}}}$}}}

\newcommand{\smcheck}[1]{\smash{\raisebox{-0.8ex}{\kern0.4ex
    \LARGE$\check{\smash{\raisebox{0.5ex}{\kern-0.2ex\normalsize
          $#1$}}}$}}}

\newcommand\subH{{\raise-1ex\hbox{$\scriptstyle H$}}}

\newcommand\barsubH{{\raise-1ex\hbox{$\left| \scriptstyle H \right|$}}}
\newcommand\barsubHo{{\raise-1ex\hbox{$\left| \scriptstyle H_0 \right|$}}}

\makeatletter
\renewcommand\section{\@startsection {section}{1}{\z@}%
                                   {-3.5ex \@plus -1ex \@minus -.2ex}%
                                   {2.3ex \@plus.2ex}%
                                   {\normalfont\normalsize\bfseries}}
\renewcommand\subsection{\@startsection{subsection}{2}{\z@}%
                                     {-3.25ex\@plus -1ex \@minus -.2ex}%
                                     {1.5ex \@plus .2ex}%
                                     {\normalfont\normalsize\bfseries}}
\renewcommand\subsubsection{\@startsection{subsubsection}{3}{\z@}%
                                     {-3.25ex\@plus -1ex \@minus -.2ex}%
                                     {1.5ex \@plus .2ex}%
                                     {\normalfont\normalsize\bfseries}}

\@addtoreset{equation}{section}
\makeatother

\newtheorem{theorem}{Theorem}[section]
\newtheorem*{theorem*}{Support theorem}

\newtheorem{corollary}[theorem]{Corollary}

\newtheorem{proposition}[theorem]{Proposition}

\renewenvironment{proof}[1][Proof]
           {\medbreak\noindent \emph{#1: \enspace}}

\newenvironment{remark}[1][Remark.]
           {\medbreak\noindent \textbf{#1 \enspace}}
           {\par \medbreak}

\makeatletter
\def\iint{\DOTSI\protect\ints@\tw@}
\def\iiint{\DOTSI\protect\ints@\thr@@}
\def\iiiint{\DOTSI\protect\ints@{4}}
\def\idotsint{\DOTSI\protect\ints@\z@}

\def\intkern@{\mkern-6mu\mathchoice{\mkern-3mu}{}{}{}}
\let\DOTSI\relax
\let\ilimits@\displaylimits

\def\ints@#1{%
  \mkern-7mu\mathchoice{\mkern-2mu}{}{}{}%
  \mathop{\mkern7mu\mathchoice{\mkern2mu}{}{}{}%
    \intop\ifnum#1=\z@\intdots@
    \else\intkern@\fi
    \ifnum#1>\tw@\intop\intkern@\fi
    \ifnum#1>\thr@@\intop\intkern@\fi
    \intop
  }\ilimits@
}
\makeatother

    {
    \begin{list}{}
        { \settowidth{\labelwidth}{}
          \leftmargin=0pt }}
    {\end{list}}

\newcommand{\romanparenlist}{
  \renewcommand{\theenumi}{\roman{enumi}}%
  \renewcommand{\labelenumi}{(\theenumi)}%
}

\newcommand{\st}[1]{\ensuremath{^{\scriptstyle \textrm{#1}}}}

\begin{document}

\begin{center}
\large{\textbf{The Abel, Fourier and Radon transforms on
    symmetric spaces}}
\end{center}
\begin{center}
{\footnotesize SIGURDUR HELGASON}
\end{center}

\begin{center}
  \textbf{To Gerrit van Dijk on the occasion of his 65\st{th} birthday}
\end{center}

\normalsize

\section{Introduction}
\label{sec:1}

In this paper we prove some recent results on the three
transforms in the title and discuss their relationships to older
results.  The spaces we deal with are symmetric spaces $X=G/K$ of
the noncompact type, $G$ being a connected noncompact semisimple
Lie group with finite center and $K$ a maximal compact subgroup.

For the two natural Radon transforms on $X$ we prove a new
inversion formula and a sharpening of an old support theorem; for
the Abel transform we prove some new identities with some
applications and for the Fourier transform a result for
integrable functions which has a strong analog of the
Riemann--Lebesgue lemma.  These latter results are from a
collaboration with Rawat, Sengupta and Sitaram.

\noindent\emph{Notation.}  Following Schwartz we use the
  notation $\D (X)$ for $\Cinf_c(X)$, $\E (X)$ for $\Cinf (X)$
  and $\SS (\RR^n)$ for the space of rapidly decreasing functions
  on $\RR^n$.

\section{Different Radon transforms on the symmetric space~$X$}
\label{sec:2}
Radon's paper \cite{R} suggested the general problem of
determining a function on a manifold on the basis of its
integrals over certain submanifolds.  A natural case of this
problem is the inversion of the X-ray transform on a Riemannian
manifold. It is the transform $f \to \hat{f}$ defined by the
arc-length integral
\begin{equation}
  \label{eq:2.1}
  \hat{f} (\gamma ) = \int_{\gamma} f(x) \, dm (x) \, ,
\end{equation}
$f$ being an ``arbitrary'' function on the Riemannian
manifold~$X$ and $\gamma$ any complete geodesic in $X$.

In general this injectivity problem seems to be unresolved.  For
a Cartan symmetric space $X \neq \SS^n$ the injectivity, however,
holds.  For a symmetric~$X$ of the noncompact the injectivity
holds in the stronger form of the 
\vspace{2ex}

\begin{theorem*} {\textup{\cite{He10}}}
\vspace{1ex}

If $\hat{f} (\gamma) =0$ for all geodesics~$\gamma$ disjoint from
a ball $B \subset X$ then $f(x) =0$ for $x \not\in B$.

\end{theorem*}

This last result requires stronger decay assumption
at~$\infty$ than the injectivity result does.

Here we shall prove an explicit inversion formula for the X-ray
transform  for
rank~$X>1$.  See \S\ref{sec:5} for the contact with Rouvi\`ere's
different solution.

Funk \cite{F} and Radon \cite{R} inverted this transform for the
sphere $\SS^2$ and $\RR^2$.  Denoting the set of geodesics by
$\Xi$ we have the coset space representations
\begin{eqnarray*}
  \SS^2 = \OO (3) /\OO (2) \, , & \Xi = \OO (3) / \OO (2) \ZZ^2\\
\RR^2 = \MM (2) /\OO (2) \, , & \Xi = \MM (2) /\MM (1) \ZZ_2
\end{eqnarray*}
$\MM (n) $ denoting the isometry group of $\RR^n$.

This suggests the following generalization.  Let $X=G/K$ and
$\Xi =G/H$ be coset spaces of the same locally compact group
$G$, $K$ and $H$ being closed subgroups.  Here it will be
convenient to assume all these groups as well as $L=K \cap H$
to be unimodular.  We do not assume the elements $\xi \in \Xi$
to be subsets of $X$ but instead use Chern's concept of incidence:
\begin{displaymath}
  x= gK \hbox{ is incident to } \xi =\gamma H
\end{displaymath}
if $gK \cap \gamma H \neq \emptyset$ as subsets of $G$.  Given
$x \in X$, $\xi \in \Xi$ define
\begin{eqnarray*}
  \bigcheck{x} &=&\{ \xi = \Xi : x\, ,\,
             \xi\hbox{ incident }\}\,  ,\\
     \hat{\xi} &=& \{ x \in X : x,\xi 
             \hbox{ incident } \} \, .
\end{eqnarray*}
These are orbits of certain subgroups of $G$ and have natural
measures $dm$, $d\mu$ (up to factors) and we define the abstract
Radon transform $f \to \hat{f}$ and its dual $\varphi \to
\checkvarphi$ by
\begin{equation}
  \label{eq:2.2}
  \hat{f} (\xi) = \int_{\hat{\xi}} f(x) \, dm (x) \, , \quad 
  \checkvarphi (x) =\int_{\bigcheck{x}} \varphi (\xi) \, 
        d \mu (\xi) \, .
\end{equation}
The normalizations of $dm$ and $d\mu$ are unified by taking $x_0
=eK$, $\xi_0 =eH$ and 
\begin{equation}
  \label{eq:2.3}
  \hat{f} (\gamma H) =\int_{H/L} f (\gamma h \cdot x_0)\, dh_L \,
  , \quad
  \checkvarphi(gK) =\int_{K/L} \varphi (gk \cdot \xi_0)
      \, dk_L
\end{equation}
the invariant measures $dh_L$, $dk_L$ being fixed by Haar
measures of $H$, $K$ and $L$.

\subsection*{Main Problems:}

\romanparenlist
\begin{enumerate}
\item 
Injectivity of $f \to \hat{f}$, $\varphi \to \checkvarphi$.

\item 
Inversion formulas.

\item 
Range and kernel question for these transforms.

\item 
Applications elsewhere.
\end{enumerate}

An easy general result relevant to problem~(iii) is the
following.  For a suitable normalization of the measures $dx=
dg_K$, $d\xi = d\gamma_H$ we have
\begin{equation}
  \label{eq:2.4}
  \int_X f(x) \checkvarphi (x) \, dx =
      \int_{\Xi} \hat{f} (\xi) \varphi (\xi) \, d\xi \, , 
\end{equation}
a result which suggests the extension of (\ref{eq:2.3}) to
distributions.

\section{$d$-planes in $\RR^n$}
\label{sec:3}

Here we consider the space $X=\RR^n$ and $\Xi = \GG (d,n)$ the
set of $d$-dimensional planes in $\RR^n$.  These are both
homogeneous under the group $G=\MM (n)$.  Fix $x_0 \in \RR^n$,
$\xi_0 \in \GG (d,n)$ at distance $d(x_0 \, , \, \xi_0)=p$.  Then
we have
\begin{equation}
  \label{eq:3.1}
  X=\RR^n =\MM (n) / K_p \, , \quad 
    \Xi = \GG (d,n) = \MM (n) / H_p\, , 
\end{equation}
where $K_p$ and $H_p$, respectively, are the stability groups of
$x_0$ and $\xi_0$.  Since various $p$ will be considered the
transforms (\ref{eq:2.2}) will be denoted by $\hat{f}_p$ and
$\sucheck{\varphi}_p$.  Since the action of $\MM (n) $ on $X$ and
$\Xi$ is quite rich it turns out that for the coset space
representation~(\ref{eq:3.1})
\begin{displaymath}
  z \in X \hbox{ is incident to } \eta \in \Xi \Leftrightarrow
  d(z,\eta) =p \, .
\end{displaymath}
Thus the transform $\hat{f}_p$ and $\checkvarphi_p$ can be
written
\begin{equation}
  \label{eq:3.2}
  \hat{f}_p (\xi) = \int_{d(x,\xi)=p} f(x) \, dm (x) \, , \, 
  \checkvarphi_p  (x) =\int_{d(x,\xi)=p}
  \varphi (\xi) \, d\mu (\xi)\, .
\end{equation}
In particular, $\hat{f}_0$ is the usual $d$-plane transform
$\hat{f}$, but in order to invert it we need $\checkvarphi_p$ for
variable $p$.  One of several versions of the inversion formula
is the following (\cite{He9}, \cite{He11}):
\begin{equation}
  \label{eq:3.3}
f(x) = c(d) \left[ \left( \frac{d}{d(r^2)} \right)^d
        \int^{\infty}_r p(p^2-r^2)^{\frac{d}{2}-1}
        (\hat{f})^{\vee}_p (x) \, dp
          \right]_{r=0}
\end{equation}
with $c(d)$ and constant.  Note that $(\hat{f})^{\vee}_p (x)$ is
the average of the integrals of $f$ over all $d$-planes at
distance $p$ from $x$.

For $d=1$ this formula reduces to
\begin{equation}
  \label{eq:3.4}
  f(x) = - \frac{1}{\pi} \int^{\infty}_0 \, 
       \frac{d}{dp} ((\hat{f})^{\vee}_p (x)) \,
          \frac{dp}{p} \, , 
\end{equation}
which for $n=2$ coincides with Radon's original formula.  Radon's
proof is very elegant and is based on an exhaustion of the
exterior $|x| >r$ by lines.  As far as I know this proof has not
been extended to higher dimensions~$n$.  Formula~(\ref{eq:3.4})
for $n>2$ is crucial for the inversion of (\ref{eq:2.1}) given in
Theorem~5.1.

\section{$d$-dimensional totally geodesic submanifolds in
  hyperbolic space $\HH^n$}
\label{sec:4}
A similar method works here and the analog of (\ref{eq:3.3}) is
the formula
\begin{equation}
  \label{eq:4.1}
  f(x) = C (d) \left[\left(  \frac{d}{d(r^2)} \right)^d 
    \int^{\infty}_r (t^2-r^2)^{\frac{d}{2}-1} t^d
      (\hat{f})^{\vee}_{s(t)} (x) \, dt \right]_{r=1}\, ,
\end{equation}
where $C (d)$ is a constant and $s(p) = \cosh^{-1} (p)$
(\cite{He9}, \cite{He11}).  Other  versions of the inversion exist
(e.g,~\cite{He12} and \cite{BC}).  For $d=1$ this reduces to
\begin{equation}
  \label{eq:4.2}
  f(x) = - \frac{1}{\pi} \int^{\infty}_0 \, 
       \frac{d}{dp} ((\hat{f})^{\vee}_p (x))\, 
           \frac{dp}{\sinh p}
\end{equation}
a formula which for $n=2$ is stated without proof in Radon \cite{R}.

\section{X-ray inversion on the symmetric space $X=G/K$}
\label{sec:5}

In communication from 2003, Rouvi\`ere proved an extension of
formula~(\ref{eq:4.2}) to symmetric spaces~$X$ of rank $\ell
=1$.  Inspired by his methods, I proved the inversion formula (\ref{eq:5.2}) for
the X-ray transform for $X$ of rank~$\ell >1$.  Then Rouvi\`ere
\cite{Ro2} extended his formula to $X$ of arbitrary rank~$\ell$.
Actually he has several such formulas but they are all different
from the formula~(\ref{eq:5.2}) below.

Fix a flat totally geodesic submanifold $E$ of $X$ with $\dim E
=\ell >1$\break ($\ell$ the rank of $X$)  passing through the origin
$o=eK$ of $X$.  Let $p>0$ and $S=S_p(o)$ be the sphere in $E$
with radius~$p$ and center~$o$.  The geodesics $\gamma$ in $E$
tangent to $S$ are permuted transitively by the orthogonal group
$\OO (E)$.  Let $du$ and $dk$ denote the normalized Haar measures
on $U$ and $K$.  The spaces $k \cdot E$ as $k$ runs through $K$
constitute all flat totally geodesic subspaces of $X$ through $o$
of dimension~$\ell$.  Thus the images $k \cdot \gamma$ ($k \in
K$, $\gamma$ tangent to $S$) constitute the set $\Gamma_p$ of all
geodesics $\gamma$ in $X$ lying in some flat $\ell$-dimensional
totally geodesic submanifold of $X$ through $o$ and $d(o,\gamma)=p$.  The set
$\Gamma_p$ has a natural measure $\omega_p$ given by the
functional
\vspace{-.5ex}
\begin{equation}
  \label{eq:5.1}
  \omega_p =\varphi \to \int_K \biggl( \int_{\OO (E)}
        \varphi (k(u \cdot \gamma)) \, du \biggr) \, dk \, .
\end{equation}

\begin{theorem}
  \label{th:5.1}
The X-ray transform (2.1) on a symmetric space
  $X=G/K$ of rank $\ell >1$ is inverted by the formula
  \begin{equation}
    \label{eq:5.2}
    f(o) =- \frac{1}{\pi} \int^{\infty}_0 \biggl( \frac{d}{dp}
         \int_{\Gamma_p} \hat{f} (\gamma) \, d\omega_p (\gamma)\biggr)
         \frac{dp}{p}\, , \quad f \in \D (X) \, .
  \end{equation}

\end{theorem}

Since $\Gamma_p$ and $d\omega_p$ are $K$-invariant the formula holds
at each point $x$ by replacing $f$ by $f \circ g$ where $g \in G$
is such that $g \cdot o=x$.

\begin{proof}
  First assume $f$ to be $K$-invariant and consider the
  restriction $f|E$.  Fix an orthonormal frame $H_0,H \in E_0$,
  the tangent space to $E$ at $o$, consider the one parameter
  subgroups $\exp tH_0$, $\exp tH$ and the geodesic $\gamma_0
  (t) = \exp tH_0 \cdot o$.  Then the geodesic $\gamma (t) =
  \exp pH \cdot \gamma_o  (t)$ lies in $E$ and is tangent to
  $S_p (o)$.  Because of (\ref{eq:3.4}) we have
  \begin{equation}
    \label{eq:5.3}
    f(o) =- \frac{1}{\pi} \int^{\infty}_0 \frac{d}{dp}\, 
    (\hat{f})^E_p (o) \frac{dp}{p}\, , 
  \end{equation}
where the superscript $E$ stands for the dual transform on
geodesics in the space~$E$.  Thus
\begin{equation}
  \label{eq:5.4}
  (\hat{f})^E_p (o) = \int\limits_{\substack{\gamma \subset E \\ d(o,\gamma)=p}}
           (\hat{f})(\gamma)\, d\nu (\gamma) =
           \int_{\OO (E)} (\hat{f}) (u \cdot \gamma)\, du \, , 
\end{equation}
where $\nu$ stands for the average over the set of
geodesics tangent to $S_p(o)$.

For $f \in \D (X)$ arbitrary we use (\ref{eq:5.2}) on the
function
\begin{displaymath}
  f^{\sharp} (x) = \int_K f(k \cdot x)\, dk \, .
\end{displaymath}
Taking into account the definition (\ref{eq:5.1}) the inversion
formula (\ref{eq:5.2}) follows immediately.

\end{proof}

\begin{remark}
  Note that the measure $\omega_p$ in (\ref{eq:5.1}) is a kind of
  convolution of the Haar measures $dk$ and $du$.  However it is
  not a strict convolution since the product $ku$ is not defined.
\end{remark}

\section{The horocycle transform in $X=G/K$}
\label{sec:6}
Consider the usual Iwasawa decomposition of $G$, $G= NAK$ where
$N$ and $A$ are nilpotent and abelian, respectively.  A
\emph{horocycle} is by definition (\cite{GG}) an orbit in $X$ of a
conjugate $gNg^{-1}$ of $N$.  The group $G$ permutes the
horocycles transitively and the space $\Xi$ of horocycles can
be written $\Xi = G/MN$ where $M$ is the centralizer of $A$ in
$K$.  In the double fibration
\begin{center}
\begin{picture}(90,50)
\put(-35,0){$X=G/K$}
\put(45,0){$G/MN =\Xi$}
\put(5,40){$G/M$}
\put(-10,15){\line(2,3){15}} 
\put(25,35){\line(2,-3){15}}
\end{picture}
\end{center}
it turns out that $x=gK$ is incident to $\xi = \gamma MN$ if and
only if $x \in \xi$.  The transforms (\ref{eq:2.2}) become
\begin{equation}
  \label{eq:6.1}
  \hat{f} (\gamma MN) =\int_N f(\gamma n \cdot o)\, dn \, , \quad 
  \checkvarphi (gK) = \int_K \varphi (gk \cdot \xi_0)\, dk\, ,
\end{equation}
where $\xi_o = N \cdot o$.  While the map $\varphi \to
\checkvarphi$ has a big kernel, the horocycle transform $f \to
\hat{f}$ is injective (\cite{He1} or \cite{GK}).  The following 
result (\cite{He5}) is considerably stronger.

\begin{theorem}{(Support theorem.)}
  \label{th:6.1}

Let $B$ be a closed ball in $X$.  Then
\begin{eqnarray*}
\begin{array}{ll}
  \hat{f} (\xi) =0 &\hbox{for } \xi \cap B =\emptyset 
       \hbox{ implies}\\
       f(x) =0 & \hbox{for } x \not\in B \, .
\end{array}
\end{eqnarray*}
\end{theorem}

Here one requires stronger decay conditions on $f$ than for the
injectivity.  A different proof was given in \cite{GQ}.  We have also the following inversion and Plancherel
formula for the Radon transform (\cite{He2,He3}).  The
pseudodifferential operator $\Lambda$ and the differential
operator~$\square$ below are constructed by means of the Harish--Chandra
$c$-function, and $w$ denotes the order of the Weyl group.  For
$G$ complex a result similar to (\ref{eq:6.2}) appears in \cite{GG}.

\begin{theorem}
  \label{th:6.2}
  \begin{trivlist}{}{}
  \item 
For $f \in \D (X)$ or sufficiently rapidly decreasing we have the
inversion formula
\begin{equation}
  \label{eq:6.2}
  f=\frac{1}{w} (\Lambda \bar{\Lambda}\hat{f})^{\vee} \, .
\end{equation}

\item 
If all Cartan subgroups of $G$ are conjugate the formula has the
improved version
\begin{displaymath}
  f=\frac{1}{w} \square ((\hat{f})^{\vee})\, .
\end{displaymath}

\item 
For $G$ arbitrary
\begin{equation}
\label{eq:6.3}
w \int_X | f(x) |^2 \, dx = \int_{\Xi} |\Lambda \hat{f}|^2
(\xi)\, d\xi\, ,
\end{equation}
with a suitable normalization of the invariant measures $dx$ and
$d\xi$.

  \end{trivlist}
\end{theorem}

The range question~(iii) for $f \to \hat{f}$ is more
complicated.  Consider first the hyperbolic plane $\HH^2$ in the
Poincar\'e unit disk model $D$.  Here the horocycles are the
circles in the disk tangential to the boundary $\{ e^{i\theta}:
\theta \in \RR \}$.  Let $\xi_{t,\theta}$ denote the horocycle
through $e^{i\theta}$ with distance $t$ (with sign) from the
origin.  Then we have the following result (\cite{He7}).

\begin{theorem}
  \label{th:6.3}
The range $\D (D)^{\hat{}}$ consists of the functions $\psi \in \D
(\Xi)$
\begin{displaymath}
  \psi (\xi_{t,\theta}) =\sum_n \psi_n (t) e^{in\theta}
\end{displaymath}
where
\begin{equation}
  \label{eq:6.4}
  \psi_n (t) = e^{-t} \biggl(\frac{d}{dt}-1\biggr) \cdots 
          \biggl(\frac{d}{dt}-2|n| +1\biggr)\varphi_n (t)
\end{equation}
where $\varphi_n \in \D (D)$ is even.
\end{theorem}

This implies a relationship between $\psi (\xi_{t,\theta})$ and
$\psi (\xi_{-t,\theta})$.  More specifically, if $f'(t) = f(-t)$,
$\Psi_n = e^t \psi_n$ then $*$ denoting convolution on $\RR$
\begin{displaymath}
  \Psi'_n = S_n * \Psi_n
\end{displaymath}
where the distribution $S_n$ on $\RR$ has Fourier transform
\begin{displaymath}
  \hat{S}_n = \frac{(i\lambda +1)\ldots (i\lambda +2|n|-1)}
        {(i\lambda -1)\ldots (i\lambda -2|n|+1)}\, , 
        \quad \lambda \in \RR\, .
\end{displaymath}
This relationship between $\psi (\xi_{-t,\theta})$ and $\psi
(\xi_{t,\theta})$ implies that in (\ref{eq:6.3}) $f \to \Lambda
\hat{f}$ does not map $L^2 (X)$ onto $L^2 (\Xi)$.

For the generalization of (\ref{eq:6.4}) to $X=G/K$ we need some
additional notation.  Let $\hat{K}$ be the unitary dual of $K$
and $d(\delta)$ the degree of a $\delta \in \hat{K}$.  Given
$\delta$ acting on $V_{\delta}$ let
\begin{displaymath}
  V^M_{\delta} = \{ v \in V_{\delta} : \delta (m) v =v
     \hbox{ for } m \in M \}
\end{displaymath}
and put $\ell (\delta)=\dim V^M_{\delta}$.  Let
\begin{displaymath}
  \hat{K}_M = \{ \delta\in \hat{K} : \ell (\delta) >0 \} \, .
\end{displaymath}
In the following theorem \cite{He9} the expansion (\ref{eq:6.5})
is a generalization of (\ref{eq:6.4}).  Put $\rho (H)
=\frac{1}{2}\Trace (\ad H |\fn )$.

\begin{theorem}
  \label{th:6.4}
The range $\D (X)\sphat$ consists of the functions $\psi \in \D
(\Xi)$
\begin{equation}
  \label{eq:6.5}
  \psi (ka \cdot \xi_0) = \sum_{\delta \in \hat{K}_M} d(\delta)
     \Tr (\delta (k) \Psi_{\delta}(a)) \qquad (\Tr = \Trace)
\end{equation}
where $\Psi_{\delta}$ is a function on $A$ with values in $\Hom
(V_{\delta}, V^M_{\delta})$, i.e.,~$\Psi_{\delta} \in \D (A,\Hom
(V_{\delta},V^M_{\delta}))$, given by
\begin{equation}
  \label{eq:6.6}
  \Psi_{\delta} (a) = e^{-\rho (\log a)} Q^{\delta}(D)_a
     (\Phi_{\delta}(a)) \qquad a \in A
\end{equation}
where
\begin{equation}
  \label{eq:6.7}
  \Phi_{\delta} \in \D (A,\Hom (V_{\delta},V^M_{\delta}))
\end{equation}
is $W$-invariant and $Q^{\delta}(D)$ is a certain $\ell (\delta )
\times \ell (\delta)$ matrix of constant coefficient differential
operators on $A$.
\end{theorem}

From this result we can derive the following (unpublished)
refinement of the support theorem above.  Let $A^+$ be the Weyl
chamber corresponding to the choice of the group~$N$.

\begin{theorem}
  \label{th:6.5}
Suppose $f \in \D (X)$ satisfies
\begin{displaymath}
\hat{f} (ka \cdot \xi_0)=0 \hbox{ for }
k \in K \, ,\,  a \in A^+\, , \, |\log a |>R \, .  
\end{displaymath}
Then 
\begin{displaymath}
  \hat{f} (ka \cdot\xi_0) =0 \hbox{ for } k \in K\, ,\, 
|\log a |>R \, ,\,  a \in A
\end{displaymath}
so by Theorem~\ref{th:6.1}
\begin{displaymath}
f(x) =0 \hbox{   for   } d (0,x)>R\, .
\end{displaymath}
\end{theorem}

\begin{proof}
Let $Q_c (D)$ be the matrix of cofactors of $Q^{\delta} (D)$ so
that
\begin{equation}
  \label{eq:6.8}
  Q_c (D) Q^{\delta}(D) = \det Q^{\delta}(D) I \, .
\end{equation}
Then (\ref{eq:6.6}) implies
\begin{equation}
  \label{eq:6.9}
  Q_c (D) (e^{\rho} \Psi_{\delta}) =\det Q^{\delta}(D)
  \Phi_{\delta}\, .
\end{equation}
Now it is known (\cite{K}, \cite{He9} pp.~267, 348) that $\det
Q^{\delta}(D)$ is a product of linear factors $\delta (H_i)+c$
where $H_i \in \fa$ and $\partial (H_i)$ the corresponding
directional derivative.

Suppose the function $\psi = \hat{f}$ satisfies
\begin{displaymath}
  \psi (ka \cdot \xi_0)=0 \hbox{   for  } k\in K \, , \, 
  a \in A^+ \, , \, |\log a| >R \, .
\end{displaymath}
Since 
\begin{displaymath}
  \Psi_{\delta}(a) = \int_K \psi (ka \cdot \xi_0)
           \delta (k^{-1}) \, dk
\end{displaymath}
we deduce from (\ref{eq:6.8}) and (\ref{eq:6.9}) that
\begin{displaymath}
  \det Q^{\delta} (D) \Phi_{\delta} (a) =0 \hbox{ for }
  a \in A^+ \, , \, |\log a |>R \, .
\end{displaymath}
Consider this equation on a ray in $A^+$ starting at $e$.
Because of the mentioned factorization of $\det Q^{\delta}(D)$ we
deduce that on this ray $\Phi_{\delta}$ satisfies an ordinary
differential equation on the interval $(R,\infty)$.  Having
compact support we deduce that $\Phi_{\delta}(a) =0$ for $a \in
A^+$, $|\log a|>R$.  By its Weyl group invariance it vanishes for
all $a \in A$, $|\log a |>R$ which by (\ref{eq:6.5}) proves the theorem.

\end{proof}

Consider the case rank $X=1$.  Let $B_R (o)$ be a ball in $X$
with radius $R$ and center $0$.  Fix a unit vector $H$  in the
Lie algebra of $A$ such that $\exp H \in A^+$.  Put $a_t=\exp
tH$.  The \emph{interior} of the horocycle $k Na_t \cdot o$ is
the union $\bigcup_{\tau >t} k N a_{\tau} \cdot o$.  A horocycle
$\xi$ is said to be \emph{external} to $B_R(o)$ if its interior
is disjoint from $B_R(o)$; $\xi$~is said to
\emph{enclose} $B_R (o)$ if its interior contains $B_R (o)$.

\begin{corollary}
  \label{cor:6.6}
Let $X$ have rank one and $B_R (o)$ as above.  Let $f \in \D
(X)$.  Then the following are equivalent:

\begin{list}{}{}
\item (i)~~$\hat{f} (\xi)=0$ whenever $\xi$ is external to
  $B_R(o)$.

\item (ii)~~$\hat{f} (\xi) = 0$ whenever $\xi$ encloses $B_R
  (o)$.

\item (iii)~$f \equiv 0$ outside $B_R (o)$.
\end{list}
\end{corollary}

For hyperbolic space this is clear from Theorem~\ref{th:6.3} and
was proved in a different way by Lax--Phillips \cite{LP}.


Problem (iii) for the dual transform $\varphi \to \checkvarphi$
has a satisfactory answer (see \cite{He9} IV \S\S 2 and 4).  The
kernel can be described in the spirit of Theorem~\ref{th:6.5} and
for the range one has the surjectivity
\begin{equation}
  \label{eq:6.10}
  \E (\Xi)^{\vee} = \E (X)\, .
\end{equation}

\section{The Abel transform}

\label{sec:7}

Let $\D_K (X)$ denote the space of $K$-invariant functions in
$\D(X)$.  The \emph{Abel transform} $f \to \A f$ is defined by
\begin{equation}
  \label{eq:7.1}
  (\A f) (a) = e^{\rho (\log a)} \int_N f(an \cdot o)\, dn\, , 
  \quad a \in A \, , f \in \D_K (X) \, .
\end{equation}
Except for the factor $e^{\rho}$ it is the restriction of the
Radon transform to $K$-invariant functions
\begin{equation}
  \label{eq:7.2}
  \A f = e^{\rho} \hat{f} \, .
\end{equation}
Some of its properties are best analyzed by means of the
spherical functions
\begin{equation}
  \label{eq:7.3}
  \varphi_{\lambda}(g) =\int_K e^{(i\lambda -\rho)(H(gk))}
      \, dk \, , \quad g \in G \, , \, \lambda \in \fa^*_c\, , 
\end{equation}
where $H(g) \in \fa$ is determined by $g \in k \exp H(g) N$ and
$\fa^*_c$ is the complex dual of $\fa$.  The \emph{spherical
  transform}
\begin{equation}
  \label{eq:7.4}
  \tilde{f} (\lambda)=\int_X f(x) \varphi_{-\lambda}(x)\, dx
  \quad f \in \D_K (X)
\end{equation}
(where $\varphi_{\lambda}(g \cdot o) = \varphi_{\lambda}(g)$) is
a homomorphism relative to convolution $\times$ on $X$:
\begin{equation}
  \label{eq:7.5}
  (f_1 \times f_2)^{\sim}(\lambda) = \tilde{f}_1 (\lambda )
        \tilde{f}_2 (\lambda )\, .
\end{equation}
As proved in \cite{H}, $\A$ intertwines the spherical transform
and the Euclidean Fourier transform $F \to F^*$ on $A$ so
\begin{equation}
  \label{eq:7.6}
  \int_A (\A f) (a) e^{-i\lambda (\log a)}\, da =
           \int_X \varphi_{-\lambda}(x)f(x)\, dx \, , 
           \quad (\A f)^* =\tilde{f} \, .
\end{equation}
Thus $\A f$ is $W$-invariant and by (\ref{eq:7.5})
\begin{equation}
  \label{eq:7.7}
  \A (f_1 \times f_2) = \A f_1 * \A f_2 \, , 
\end{equation}
where $*$ is convolution on $A$.  Let $\DD (X)$ denote the
algebra of $G$-invariant differential operators on $X$ and
$\Gamma : \DD (X) \to \DD_W (A)$ the isomorphism onto the
$W$-invariant constant coefficient differential operators on~$A$.

The Abel transform is a simultaneous transmution operator between
$\DD (X)$ and $\DD_W(A)$, i.e.,
\begin{equation}
  \label{eq:7.8}
  \A Df = \Gamma (D) \A f \, , \quad D \in \DD (X) \, , \, 
  f \in \D_K (X)
\end{equation}
(\cite{He13}) which for example can be used to prove that each
$D$ has a fundamental solution.  By the Paley--Wiener theorem for
(\ref{eq:7.3}) one has that $\A :\D_K(X) \to \D_W (A)$ is a
bijective homeomorphism.  (Here the subscript $W$ means
$W$-invariance.)  Hence we have a bijection
\begin{equation}
  \label{eq:7.9}
  \A^* : \D'_W (A) \to \D'_K (X)
\end{equation}
between the corresponding distribution spaces.  Also if $\varphi
\in \E_W (A)$ we have easily (\cite{Be} or \cite{He9} IV, \S4)
\begin{equation}
  \label{eq:7.10}
  (\A^* \varphi) (gK) = \int_{K/M} \varphi (\exp H(gk))
      e^{-\rho (H(gk))} \, dk \, .
\end{equation}
The Radon transform has the advantage over $\A$ that it commutes
with the action of $G$.  Thus we can deduce from (\ref{eq:6.10})
and (\ref{eq:7.2}) that (\cite{He9})
\begin{equation*}
  \A^* \E_W (A) = \E_K (X) \, .
\end{equation*}
We now add a few new results about $\A$ and $\A^*$ which will be
useful later.  Some are closely related to rank-one results in
Bagchi and Sitaram in \cite{BS}.

Because of the convolution property (\ref{eq:7.7}) one can ask
how $\A^*$ behaves relative to convolution.  Let $L$ be the
operator on $\s (A)$ given by
\begin{equation}
  \label{eq:7.11}
  (L \varphi)^* (\lambda) = |c(\lambda)|^{-2} \varphi^*
      (\lambda), \quad \lambda \in \fa^* \, ,
\end{equation}
where $c(\lambda)$ is Harish--Chandra's $c$-function.

\begin{theorem}
  \label{th:7.1}
Let $\varphi \in \D_W (A)$, $\psi \in \E_W (A)$.  Then
\begin{displaymath}
  \A^* (L\varphi) =w \,\, \A^{-1}(\varphi) \qquad
     (w = \hbox{ order of }W)
\end{displaymath}
and
\begin{equation}
\label{eq:7.12}
  \A^* (\varphi * \psi) = \frac{1}{w}\A^* (L\varphi)
      \times \A^* \psi \, .
\end{equation}

\end{theorem}

\begin{proof}
Using the inversion formula for the spherical transform we have
\begin{eqnarray*}
  \A^* (L\varphi)(gK) &=& \int_K (L\varphi)
     (\exp H(gk))^{-\rho (H(gk))}\, dk\\
     &=& \int_K \biggl(\int_{\fa^*} (L\varphi)^* (\lambda)
          e^{i\lambda (H(gk))}\, d\lambda \biggr) 
          e^{-\rho (H(gk))}\, dk\\
     &=& \int_{\fa^*} |c(\lambda)|^{-2} \varphi^*_{\lambda}
          (g) \, d\lambda =F(gK)\, , \\
\end{eqnarray*}
where $\tilde{F} (\lambda) = \varphi^* (\lambda) w$.  But 
$\tilde{F} = (\A F)^* =\varphi^* w$ so
\begin{displaymath}
  \varphi = \frac{1}{w}\A F \, , \quad F=w \A^{-1}(\varphi)\, .
\end{displaymath}
Thus  $\A^* (L\varphi) = w \A^{-1}(\varphi)$.  Consider now the average
\begin{displaymath}
  \psi^{\lambda}(a) = \frac{1}{w}\sum_{s \in W}
      e^{is\lambda (\log a)}\, .
\end{displaymath}
Then $  \A^* \psi^{\lambda}=\varphi_{\lambda}$ and  
  $\varphi * \psi^{\lambda}=\varphi^* (\lambda)\psi^{\lambda}$.  But $\A^{-1}\varphi = \frac{1}{w} F \in \D_K(X)$ and
\begin{displaymath}
  \A^{-1} \varphi \times \varphi_{\lambda} =\frac{1}{w}
     \tilde{F} (\lambda) \varphi_{\lambda} =\varphi^*
        (\lambda) \varphi_{\lambda}\, .
\end{displaymath}
Combining these formulas we have
\begin{equation}
  \label{eq:7.13}
  \A^* (\varphi * \psi^{\lambda}) = \A^{-1}\varphi
      \times \A^* (\psi^{\lambda})\, .
\end{equation}

Now $\psi \in \D_W (A)$ is a superposition
\begin{displaymath}
  \psi (a) =\int_{\fa^*}\psi^* (\lambda)\psi^{\lambda}
      (a) \, d\lambda
\end{displaymath}
so the identity (\ref{eq:7.12}) follows from (\ref{eq:7.13}) for
such $\psi$. For $\psi \in \E_W(A)$ the identity follows by an
approximation because $\varphi$ and $\A^* (L\varphi)$ have
compact support and $\A^*$ is continuous on $\E_W(A)$. 

\end{proof}

Theorem~\ref{th:7.1} implies the following inversion formula
which in reality is a special case of (\ref{eq:6.2}).   It
appears also in \cite{Be}.

\begin{corollary}
  \label{cor:7.2}
The transform $f \to \A f$ has inversion
\begin{displaymath}
  f=\frac{1}{w} \A^* (L\A f)\, \qquad f \in \D_K(X)\, .
\end{displaymath}

\end{corollary}

The above results suggest various ways of defining $\A$ on the
space $\E'_K(X)$ of $K$-invariant compactly supported
distributions on $X$ although formula (\ref{eq:7.1}) does not
work.

\begin{trivlist}{}{}
\item 
\emph{Spherical transform method.}\quad If $T \in \E'_K(X)$, the
spherical transform
\begin{displaymath}
  \tilde{T} (\lambda) = \int_X \varphi_{-\lambda} (x)
      \, dT (x)
\end{displaymath}
is a $W$-invariant entire function of exponential type on
$\fa^*_c$ and of polynomial growth.  (See \cite{EHO} or
\cite{He5}, Theorem~8.5.)  By the Euclidean
Paley--Wiener theorem there exists an $S \in \E'_W (A)$ such that
$\tilde{T} =S^*$.  Thus in accordance with (\ref{eq:7.6}) we put
\begin{equation}
  \label{eq:7.14}
  \A T=S \, .
\end{equation}

\item 

\emph{Radon transform method.}\quad Because of (\ref{eq:2.4}) the
Radon transform of a distribution $T \in \E'(X)$ is defined by
\begin{displaymath}
  \hat{T} (\varphi) = T(\checkvarphi)\qquad \varphi \in \E
  (\Xi)\, .
\end{displaymath}
If $T$ is $K$-invariant then so is $\hat{T}$ and since $\Xi
=K/M \times A$ under the bijection $(kM,a) \to ka \cdot
\xi_o$ we see that $\hat{T}$ has the form $\hat{T} = 1
\otimes \sigma$ where $\sigma \in \E'(A)$.  Because of (\ref{eq:7.2}) we put
\begin{equation}
  \label{eq:7.15}
  \A T = e^{\rho} \sigma \, .
\end{equation}

\item 
\emph{Functional analysis method.}\quad As remarked $\A^*$ is a
bijection of $\D'_W (A)$ onto $\D'_K(X)$.  The restriction of
$\A^*$ to $\E_W(A)$ is a continuous bijection onto $\E_K(X)$ and
in fact a homeomorphism since both spaces are Fr\'echet.  Thus we
have $(\A^*)^* : \E'_K(X) \to \E'_W (A)$ bijectively so we can
define
\begin{equation}
  \label{eq:7.16}
  \A T = (\A^*)^* (T) \, .
\end{equation}

\end{trivlist}

\begin{proposition}
  \label{prop:7.3}
All the definitions (\ref{eq:7.14})--(\ref{eq:7.16}) coincide.
\end{proposition}

The convolution property in Theorem~\ref{th:7.1} extends readily
to distributions so
\begin{eqnarray*}
  \A^* (\E'_W (A)*\psi) &=& \A^{-1}(\E'_W (A))
       \times \A^* \psi\\
       &=& \E'_K(X) \times \A^* \psi \, .
\end{eqnarray*}
Thus putting
\begin{eqnarray*}
  V_{\psi} = \E'_W(A) * \varphi \, , \quad
             \psi \in \E_W (A)\, ; \quad 
  W_f = \E'_K(X) \times f \, , \quad
        f \in \E_K (X)
\end{eqnarray*}
we conclude that
\begin{equation}
  \label{eq:7.17}
\A^* (V_{\psi}) = W_{\A^* \psi}  \, .
\end{equation}

\begin{theorem}[(Bagchi--Sitaram)]
  \label{th:7.4}
If $X$ has rank one and $f \in \E_K(X)$ then the closure of the
space $W_f =\E'_K(X) \times f$ contains a spherical function.

\end{theorem}

The authors use (\ref{eq:7.17}) to reduce the question to the
analogous one for the one-dimensional space $A \sim \RR$ where by
Schwartz's theorem stated in \S9 below some exponentials $e^{i\mu}$ and $e^{-i\mu}$
belong to the closure and $\A^* (e^{i\mu} + e^{-i\mu}) = 2\varphi_{\mu}$.

\section{The Fourier transform on $X=G/K$}
\label{sec:8}

We now go to the notation of \S6 with the Iwasawa decomposition
$G= NAK$, and $\fg =\fn +\fa +\fk$ for the corresponding Lie
algebras. For $g \in G$ let $A(g) \in \fn$ be determined by $g=n
\exp A(g)k$ $(n \in N , k \in K)$.  Given $x=gK$ in $X$, $b=kM$
in $B=K/M$ we put
\begin{displaymath}
  A(x,b) = A (k^{-1}g)
\end{displaymath}
and as usual we put $\rho (H) =\frac{1}{2}\Trace (\ad H|\fn)$.
Let $\fa^*_c$ denote the space of complex-valued linear forms on
$\fa$.

Given a function $f$ on $X$ we define its \emph{Fourier
  transform} by
\begin{equation}
  \label{eq:8.1}
  \tilde{f} (\lambda ,b) = \int_X f(x)
      e^{(-i\lambda +\rho) (A (x,b))}\, dx
\end{equation}
for those $(\lambda ,b) \in \fa^*_c \times B$ for which the
integral is defined.  Many of the principal theorems for Fourier
transforms on $\RR^n$ have analogs for $X=G/K$.

\begin{trivlist}{}{}
\item \emph{Inversion Formula}(\cite{He3}).\quad For $f \in
  \D(X)$ we have
  \begin{displaymath}
    f(x) =\frac{1}{w}\int_{\fa^*} \!\! \int_{B}\tilde{f}
       (\lambda,b)e^{(i\lambda +\rho) (A(x,b))}
       |c(\lambda)|^{-2}\, d\lambda \, db
  \end{displaymath}
where $c(\lambda)$ is Harish--Chandra's $c$-function.

\item  \emph{Plancherel Formula}(\cite{He4}).\quad The map $f \to
  \tilde{f}$ extends to an isometry of $L^2(X)$ onto $L^2
  (\fa^*_+ \times B)$:
  \begin{equation}
    \label{eq:8.2}
    \int_X |f(x)|^2 \, dx = \int_{\fa^*_+ \times B}
    |\tilde{f} (\lambda ,b|^2)|c(\lambda)|^{-2}\, 
    d\lambda \, db \, .
  \end{equation}

\item   \emph{Paley-Wiener Theorem}(\cite{He5}).\quad The map $f
  \to \tilde{f}$ maps the space $\D(X)$ onto the space of smooth
  $\varphi (\lambda ,b)$ on  $\fa^*_c \times B$ which are
  holomorphic on $\fa^*_c$ of exponential type (uniformly in $B$)
  satisfying the invariance condition
  \begin{equation}
    \label{eq:8.3}
    \int_B \varphi (\lambda ,b) 
       e^{(i\lambda +\rho)(A(x,b))}\, db \,\,
\hbox{ is $W$-invariant in } \lambda \, .
  \end{equation}

\end{trivlist}

For the next result we refer to Eguchi's paper for full
explanations of notation.

\begin{trivlist}{}{}
\item \emph{The Schwartz Theorem}(\cite{Eg}).\quad Let $0<p\leq2$
  and $\s^p(X) \subset L^p(X)$  the corresponding Schwartz
  space.  Let $\epsilon = \frac{2}{p}-1$ and $\s (\fa^*_{\epsilon}
  \times B)$ the space of functions which are holomorphic in the
  ``tube'' $\fa^*_{\epsilon} \times B$, are rapidly decreasing
  and satisfy (\ref{eq:8.3}).  Then $f \to \tilde{f}$ is a
  bijection of $\s^p(X)$ onto $\s(\fa^*_{\epsilon} \times B)$.

\end{trivlist}

These results leave out the space $L^1 (X)$ and one should think
that a self-respecting Fourier transform should be defined here.

We shall now show (modifying a bit the proof of \cite{HRSS}) that this
can be done and that a strong analog of the classical
Riemann--Lebesgue lemma holds for $\tilde{f}$ in (\ref{eq:8.1}).

Let $C(\rho)$ denote the convex hull in $\fa^*$ of the set $\{ s
\rho : s \in W \}$ of Weyl group transforms of $\rho$.

\begin{theorem}
  \label{th:8.1}
Let $f \in L^1 (B)$.  Then there exists a subset $B' \subset B$
with $B-B'$ of measure $0$ such that for each $b \in B'$

\begin{trivlist}{}{}
\item 
(i)~~$\tilde{f} (\lambda ,b)$ is defined for $\lambda$ in the tube
$\fa^* +i C(\rho)$ and holomorphic in its interior.

\item 
  (ii)~~$\lim_{\xi \to \infty}\, \tilde{f} (\xi + i \eta ,b)=0$
  uniformly for $\eta \in C(\rho)$.
\end{trivlist}
\end{theorem}

\begin{proof}
Let $\lambda =\xi +i\eta$ where $\xi \in \fa^* \, , \, \eta \in C
(\rho)$.  Then 
\begin{equation}
  \label{eq:8.4}
  \int_B |\tilde{f} (\lambda ,b)| \, db \leq \int_X |f(x)|
      \int_B e^{(\eta + \rho)(A(x,b))}\, db \, dx \, .
\end{equation}
The integral over $B$ is the spherical function
$\varphi_{-i\eta}$ which is bounded by~$1$ (\cite{HJ}).  Thus
\begin{displaymath}
  \| \tilde{f} (\lambda , \cdot) \|_1 \leq \| f \|_1
\end{displaymath}
and for each $\lambda \in \fa^* +i C(\rho)$,
$\tilde{f} (\lambda ,b)$ exists for all $b$ in a subset
$B_{\lambda} \subset B$ of full invariant measure.  Let
\begin{displaymath}
  B' = B' (f) =\cap_{s \in W} B_{is\rho}\, .
\end{displaymath}
For the statements~(i) and (ii) we may assume $f \geq 0$ in
(\ref{eq:8.1}).  Since $b \in B_{is\rho}$ for each $s \in W$ we
have
\begin{equation}
  \label{eq:8.5}
  \int_X f(x) e^{(s\rho +\rho)(A(x,b))} \, dx <\infty \, .
\end{equation}
Fix $b \in B'$, $\eta \in C(\rho)$.  Then 
\begin{equation}
  \label{eq:8.6}
  \int_X f(x) e^{(\rho + \eta) (A(x,b))} \, dx =
     \sum_{\sigma \in W} \int_{X_{\sigma}} f(x)
       e^{(\rho + \eta)(A(x,b))} \, dx \, 
\end{equation}
where 
\begin{displaymath}
  X_{\sigma}=\{ x \in X : \sigma ( A(x,b)) \in \overline{\fa^+}\,
  .
\end{displaymath}
Replace $\eta (A(x,b))$ by $(\sigma \eta) $ $(\sigma(A(x,b)))$ and
let $(\sigma \eta)^+$ be the element in $\overline{\fa^*_+}$,
which is $W$-conjugate to $\sigma \eta$.  Then since $(\sigma
\eta)^+-\sigma \eta \geq 0$ on $\fa^+$ we have
\begin{displaymath}
  \int_{X_{\sigma}} f(x) e^{(\rho +\eta)(A(x,b))} \, dx \leq
    \int_{X_{\sigma}} f(x) 
    e^{(\rho + (\sigma \eta)^+)(\sigma (A(x,b)))}\, dx \, .
\end{displaymath}
Now by Lemma~8.3, Ch.~IV in \cite{He8}
\begin{displaymath}
  \overline{\fa^*_+} \cap C(\rho) = 
     \overline{\fa^*_+} \cap (\rho + {}^- \fa^*) \, ,
\end{displaymath}
where
\begin{displaymath}
  {}^-\fa^* = \{ \lambda \in \fa^* | \langle \lambda ,\mu \rangle
     \leq 0 \hbox{ for } \mu \in \fa^*_+ \} \, .
\end{displaymath}
Thus
\begin{displaymath}
  (\sigma \eta)^+ \in \overline{\fa^*_+} \cap
  (\rho + {}^-\fa^*)\,,
\end{displaymath}
whence
\begin{equation}
  \label{eq:8.7}
  (\sigma \eta)^+ -\rho \leq 0 \hbox{ on } \fa^+ \, .
\end{equation}
Thus the last integral is bounded by
\begin{displaymath}
  \int_{X_{\sigma}} f(x) 
     e^{(\rho +\sigma^{-1}\rho)(A(x,b))} \, dx <\infty
\end{displaymath}
by (\ref{eq:8.5}).  This shows by (\ref{eq:8.6}) that if $b \in
B'$ and $\lambda \in \fa^* +i C (\rho)$ the integral
(\ref{eq:8.1}) is absolutely convergent.  The holomorphy
statement follows by Morera's theorem.  This proves~(i).

For part (ii) we use the Radon transform (\ref{eq:6.1}).  Since
$f \in L^1(X)$, $\hat{f} (\xi)$ exist for almost all $\xi \in
\Xi$ (\cite{He4}).  Since $(kM,a) \to ka \cdot \xi_0$ is a
diffeomorphism of $K/M$ onto $\Xi$ we write $\hat{f}(kM,a)$ for
$\hat{f} (ka \cdot \xi_0)$.  Enlarging $B'$ to another subset of
$B$ of full invariant measure we may assume $\hat{f} (b,a)$
exists for $b \in B'$ and almost all $a$.  Now we have
\begin{equation}
  \label{eq:8.8}
  \int_X f(x) \, dx = \int_{AN} f(anK)\, da \, dn
\end{equation}
for suitable Haar measures on $A$ and $N$.  Applying this to the
function $x \to f(k \cdot x)$ with $kM =b \in B'$ we get
\begin{equation}
  \label{eq:8.9}
  \int_X f(x) \, dx =\int_A \hat{f} (kM,a)\, da
\end{equation}
so since $A (an \cdot o)= \log a$,
\begin{eqnarray}
  \label{eq:8.10}
  \tilde{f} (\lambda ,kM) &=& \int_A \hat{f} (kM,a)
     e^{(\rho +\eta) (\log a)} e^{-i\xi (\log a)}\, da\\
\nonumber
     &=& \sum_{s \in W} \int_{s^{-1}A^+} \hat{f}
         (kM,a) e^{(\rho +\eta)(\log a)}
         e^{-i\xi(\log a)}\, da\, .
\end{eqnarray}

Now $a \in s^{-1} A^+$ implies $sa \in A^+$ and 
\begin{displaymath}
  \eta (\log a) = (s\eta) (s\log a) \leq (s\eta)^+
      (s\log a) \leq \rho (s \log a)
\end{displaymath}
by (\ref{eq:8.7}).  Thus on $s^{-1} A^+$, 
\begin{equation}
  \label{eq:8.11}
  \hat{f} (kM,a) e^{(\rho + \eta)(\log a)} \leq
  \hat{f} (kM,a) e^{(\rho + s^{-1}\rho)(\log a)}\, .
\end{equation}
For $b=kM \in B'$ the integral in (\ref{eq:8.1}) is absolutely
convergent so by (\ref{eq:8.9}) the function
\begin{equation}
  \label{eq:8.12}
  a \to \hat{f} (kM,a) e^{(\rho + \eta)(\log a)}
\end{equation}
belongs to $L^1 (A)$.  The first part of (\ref{eq:8.10}) combined
with the Riemann--Lebesgue lemma for the Fourier transform on $A$
shows that for each $\eta \in C (\rho)$
\begin{displaymath}
  \lim_{\xi \to \infty} \tilde{f} (\lambda ,b)=0\, .
\end{displaymath}
For the uniform convergence in (ii) we use the second part of
(\ref{eq:8.10}).  Let $f_n$ be positive in $\D (X)$ such that $f_n
\to f$ a.e. and $f_n (x) \leq f(x)$.  In (\ref{eq:8.10}) and
(\ref{eq:8.11}) we replace $f$ by the function $g_n =f-f_n$.
Then 
\begin{eqnarray*}
  |\tilde{g}_n (\lambda ,kM) | & \leqq & \sum_{s \in W}
     \int_{s^{-1}A^+} \hat{g}_n (kM,a)
        e^{(\rho +\eta)(\log a)}\, da \\
        &\leq & \sum_{s \in W} \int_{s^{-1}A^+}
          \hat{g}_n (kM,a) e^{(\rho +s^{-1}\rho)(\log a)}\, 
             da \, ,
\end{eqnarray*}
which tends to $0$ as $n \to \infty$ by (\ref{eq:8.9}),
(\ref{eq:8.12})  and the dominated convergence theorem.  Thus
given $\epsilon >0$ we can fix $N$ such that $|\tilde{g}_N
(\lambda ,kM)| <\epsilon$ for all $\lambda \in \fa^* +iC_{\rho}$.  By
the Paley--Wiener theorem for $\D (X)$ there is an $L$ such that
$|\tilde{f}_N (\xi +i\eta ,kM)|\leq \epsilon$ for $|\xi | >L$ and
$\eta \in C (\rho)$.  Since $\tilde{g}_N = \tilde{f}-\tilde{f}_N$
 this proves (ii).

\end{proof}

\begin{remark}
  Another version of (ii) involving the $L^1$ norm over $B$ is
  given in \cite{SS}.
\end{remark}

\section{Spectral analysis on $X$}
\label{sec:9}

A theorem of Schwartz \cite{Sc} states that if $f$ is a function
in $\E (\RR)$ $(f \not\equiv 0)$ the closed subspace of $\E(X)$ (in
its usual Fr\'echet space topology) generated by all the
translates of $f$ contains an exponential $e^{\mu x}$ for some
$\mu \in \CC$.

We shall now give the proof from \cite{HS} of the following
analog of Schwartz's theorem.

\begin{theorem}
  \label{th:9.1}
Let $X =G/K$ have rank one and $f\neq 0$ a function in $\E(X)$.
Then the closed subspace $V_f$ of $\E(X)$ generated by the
$G$-translates of $f$ contains a function
\begin{displaymath}
  x \to e_{\mu ,b} (x) =  e^{\mu (A(x,b))}
\end{displaymath}
for some $\mu \in \fa^*_c$.
\end{theorem}

For this we consider for $\lambda \in \fa^*_c$ the \emph{Poisson
transform}
\begin{displaymath}
  \P_{\lambda}: F(b) \to f(x)\, , \, \qquad F \in L^1(B)\, , 
\end{displaymath}
where
\begin{equation}
  \label{eq:9.1}
  f(x) =\int_B e^{(i\lambda +\rho)(A(x,b))}F(b) \, db \, .
\end{equation}
The element $\lambda$ is said to be \emph{simple} if $\P_{\lambda}$
is injective.  The simplicity criterion for $\lambda$ \cite{He6}
implies that for each $\lambda \in \fa^*_c$ one of the transforms
$s\lambda$\,\,\, $(s \in W)$ is simple.  Consider now the spherical
function $\varphi_{\lambda}$  (\ref{eq:7.3}) which can also be
written
\begin{displaymath}
  \varphi_{\lambda} (x) =\int_B 
       e^{(i\lambda +\rho)(A(x,b))}\, db \, .
\end{displaymath}
We know from \cite{He9}, III, Lemma~2.3 that if $-\lambda$ is
simple then the closed space $\E_{(\lambda)}(X) \subset \E (X)$
generated by the $G$-translates of $\varphi_{\lambda}$ contains
the space $\P_{\lambda} (L^2 (B))$.

Coming to the proof of the theorem we conclude from the
Bagchi--Sitaram result (\ref{eq:7.4}) that the space
$V^K_f$ of $K$-invariants in $V_f$ contains a spherical
function $\varphi_{\lambda}$.  By the simplicity result quoted,
either $\lambda$ or $-\lambda$ is simple so we can take
$-\lambda$ simple.  Thus by the conclusion above, $V_f$ contains
the space $\P_{\lambda} (L^2 (B))$.  Now by \cite{He9}, III,
Exercise~B1, pp.~371 and 570,
\begin{equation}
  \label{eq:9.2}
  e^{(i\lambda +\rho) A(x,eM)} =\sum_{\delta \in \hat{K}_M}
     \, d(\delta) \varphi_{\lambda,\delta}(x) \, , 
\end{equation}
with $\delta$ and $\hat{K}_M$ as in (\ref{eq:6.5}) and
\begin{displaymath}
  \varphi_{\lambda ,\delta} (x) = \int_{K}
  e^{(i\lambda +\rho) (A (x,k))}\langle \delta (k) v,v \rangle
   \, dk \, .
\end{displaymath}
Thus $\varphi_{\lambda ,\delta} \in V_f$ so since (\ref{eq:9.2})
converges in the topology  of $ \E (X)$ the theorem follows.

\begin{remark}
  Since Schwartz's theorem fails for $\RR^n$ $(n>1)$ the proof
  above via the Bagchi--Sitaram theorem is limited to the case of
  rank~$X=1$.  However, this does not rule out the possibility
  that Theorem~\ref{th:9.1} might remain valid for $X$ of higher rank.
\end{remark}

\section{Further results on the Fourier transform}
\label{sec:10}

A result of Hardy's \cite{Ha} shows limitations on how fast a
function on $\RR^n$ and its Fourier transform can decay at
$\infty$.  Precisely, if
\begin{displaymath}
  |f(x)|\leq A e^{-\alpha |x|^2}\, , \, 
  |\tilde{f} (u) |\leq B e^{-\beta |u|^2}\quad
  \alpha ,\beta >0
\end{displaymath}
and if $\alpha \beta >\frac{1}{4}$ then $f=0$.  Sitaram and
Sundari \cite{SiSu} proved an analog for a class of spaces $X$
and Sengupta \cite{Se} extended this to all~$X$.  Many other
variations of the result have been proved by Ray and Sarkar,
Cowling, Sitaram and Sundari, Narayanan and Ray, Shimeno,
Thangavelu.  (See References.)

The following classical result is closely related to Wiener's
Tauberian theorem.  \emph{Let $f \in L^1 (\RR^n)$ such that
$\tilde{f} (u) \neq 0$ for all $u \in \RR^n$.  Then the
translates of $f$ span a dense subspace of $\RR^n$}.  Many papers
deal with analogies of this result for semisimple Lie groups and
symmetric spaces.  See \cite{EM}, \cite{Sa}, \cite{Si1},
\cite{Si2}, \cite{SS}, \cite{MRSS} for a sample.

The polar coordinate representation $(kM,a)\to kaK$ of $X$
identifies $X$ with $K/M \times A^+$ up to a null set.  Thus one
might interpret the Plancherel formula (\ref{eq:8.2}) as
identifying $X$ with its ``dual''.  But in contrast to $\RR^n$
where the Fourier transform is essentially equal to its inverse,
the Fourier transform
\begin{eqnarray}
  \label{eq:10.1}
  \tilde{f} (\lambda ,b) &=&  \int_{X} f(x)
     e^{(-i\lambda +\rho)(A(x,b))}\, dx \, ,\\
\noalign{\nonumber\hbox{and the inverse}}\\
\label{eq:10.2}
(\F^{-1}\varphi) (x) &=& \int_{\fa^* \times B}
   \varphi (\lambda ,b) e^{(i\lambda +\rho)(A(x,b))}
   |c(\lambda)|^2 \, d \lambda \, db
\end{eqnarray}
are quite different.  Hence it is a natural problem to prove the
analog of the Paley--Wiener theorem for $\F^{-1}$.

This was done by A.~Pasquale \cite{Pa} for the spherical
transform for $X$ of rank one or the case of $G$ complex, and by
N.~Andersen \cite{AN} in general.  Let $L$ denote the Laplacian
on $X$.

\begin{theorem}
  \label{th:10.1}
The image of $\F^{-1} (\D (\fa^* \times B))$ consists of the
functions $f$ on $X$ satisfying
\begin{displaymath}
  (1+ \, d(o,x))^m L^n f \in L^2 (X) \quad
  \hbox{for all } m,n \in \ZZ^+
\end{displaymath}
and
\begin{displaymath}
  \lim_{n \to \infty} \| (L+\langle\rho ,\rho\rangle )^n
       \|^{1/2n}_2 < \infty \, .
\end{displaymath}

\end{theorem}

Another characterization was given by Pesenson \cite{P}, namely
\begin{displaymath}
  \| L^{\sigma}f\|_2 \leq (\omega^2 +|\rho |^2)^{\sigma}
  \| f \|_2 \quad \hbox{for all } \sigma >0 \, .
\end{displaymath}

\end{document}